\RequirePackage{ifpdf}
\ifpdf 
\documentclass[pdftex]{sigma}
\else
\documentclass{sigma}
\fi

\newcommand{\C}{\ensuremath{\mathbb C}}

\newcommand{\Z}{\ensuremath{\mathbb Z}}

\begin{document}

\allowdisplaybreaks
	
\renewcommand{\PaperNumber}{101}

\FirstPageHeading

\renewcommand{\thefootnote}{$\star$}

\ShortArticleName{Toeplitz Operators, K\"ahler
Manifolds, and Line Bundles}

\ArticleName{Toeplitz Operators, K\"ahler
Manifolds,\\ and Line Bundles\footnote{This paper is a
contribution to the Proceedings of the 2007 Midwest
Geometry Conference in honor of Thomas~P.\ Branson. The full collection is available at
\href{http://www.emis.de/journals/SIGMA/MGC2007.html}{http://www.emis.de/journals/SIGMA/MGC2007.html}}}

\Author{Tatyana FOTH}

\AuthorNameForHeading{T.~Foth}

\Address{Department of Mathematics,
University of Western Ontario,\\
London, Ontario N6A 5B7, Canada}

\Email{\href{mailto:tfoth@uwo.ca}{tfoth@uwo.ca}}

\URLaddress{\url{http://www.math.uwo.ca/~tfoth/}}

\ArticleDates{Received August 23, 2007, in f\/inal form October 23, 2007; Published online October 26, 2007}

\Abstract{This is a survey paper. We discuss
Toeplitz operators in K\"ahler
geometry, with applications to geometric quantization,
and review some recent developments.}


\Classification{32Q15; 53D50}

\section{Historical remarks}

Toeplitz operators have been studied by analysts
for many years.
In \cite{brown:63} a {\it Toeplitz operator}
on $S^1=\{ z\in\C | \ |z|=1\}$ is def\/ined as follows.
Denote by $\mu$ the normalized Lebesgue measure,
and by $e_n=e_n(z)=z^n$, $z\in S^1$, $n\in\Z$.
Functions $e_n$ are bounded, measurable,
and they form an orthonormal basis in $L^2=L^2(S^1,\mu)$.
A function $f\in L^2$ is called {\it analytic}
if $\int_{S^1} f\bar{e}_nd\mu=0$ for all $n<0$,
and the {\it Hardy space}
$H^2$ is def\/ined
as the space of all functions in $L^2$ which are
analytic. Denote by $P :L^2\to H^2$ the orthogonal
projector. Let $\varphi$ be a bounded measurable function
on~$S^1$. The corresponding {\it Toeplitz operator}
$T_{\varphi}:H^2\to H^2$ is def\/ined by
$T_{\varphi}=P\circ M_{\varphi}$, where $M_{\varphi}$
is the operator of multiplication by $\varphi$.
The function $\varphi$ is called {\it the symbol of}
$T_{\varphi}$. We immediately observe:
for $\varphi(z)=1$ $T_{\varphi}$
is the identity operator, and for
$\alpha, \beta\in \C$ and bounded measurable functions
$f,g$ on $S^1$ we have:
$T_{\alpha f+\beta g}=\alpha T_f+\beta T_g$.
These operators have various remarkable properties.
For example, Theorem $4$
\cite{brown:63} states:
A necessary and suf\/f\/icient condition
that an operator on $H^2$ be a Toeplitz operator
is that its matrix (with respect to the orthonormal
basis $\{ e_n: \ n=0,1,2,\dots\}$) be a Toeplitz matrix.

\begin{remark} Otto Toeplitz (1881--1940)
was a German born Jewish mathematician, professor in Bonn
from 1928 until 1935. A {\it Toeplitz matrix}
is a one-way inf\/inite matrix $(a_{ij})$
(i.e. $i,j=0,1,2,\dots$) such that $a_{i+1,j+1}=a_{ij}$.
For example, the matrix of $T_{e_1}$
in the basis $\{ e_n: \ n=0,1,2,\dots\}$ is
\[
\begin{pmatrix}
0 & 0 & 0 & 0 & \cdots \\
1 & 0 & 0 & 0 & \cdots \\
0 & 1 & 0 & 0 & \cdots \\
0 & 0 & 1 & 0 & \cdots \\
\cdots & \cdots & \cdots & \cdots & \cdots
\end{pmatrix}.
\]
\end{remark}

Toeplitz operators on
bounded domains
in $\C^n$, $n\ge 1$, (on various function spaces)
have been studied extensively,
and it would be a very dif\/f\/icult task to give a
comprehensive description of all the work done
in this area. See, in particular,
\cite{bottcher:90,douglas:05,douglas:70,jewell:80,jewell:79,
klein:72,klein:73,peller:93,salinas:89,stroethoff:03, stroethoff:07,upmeier:83,upmeier:84,upmeier:85,upmeier:96,
widom:64,widom:66},
and references in \cite{bottcher:90}.

The purpose of this paper is dif\/ferent:
we outline how Toeplitz operators appear in complex
and symplectic geometry (often in problems
coming from mathematical physics)
and we overview some recent results in this area.

\section{K\"ahler manifolds, geometric quantization\\
and Toeplitz operators}

\subsection{Preliminaries}
\label{bg}
The use of Toeplitz operators in geometric quantization
has its origins in work
of F.~Berezin, L.~Boutet de Monvel, and
J.~Sj\"ostrand, see, in particular,
\cite{berezin:72,berezin:74,berezin:75,
boutet:78,boutet:81,boutet:76},
and see also~\cite{kostant:70,souriau:97}.

Many of the key ideas are contained in \cite{boutet:78}.
The article \cite{boutet:76} laid down
the foundations for the analysis.
Other important, more recent, papers include
\cite{berman:05} and \cite{boutet:05}.

In the following ``smooth'' will always mean $C^{\infty}$.
Let $W$ be a strictly pseudoconvex domain in
a complex $n$-dimensional manifold, $n\ge 1$.
Assume that the boundary $\partial W$ is smooth
and $\overline{W}=W\cup \partial W$ is compact.
Let $r\in C^{\infty} (M)$ be a def\/ining function
for $W$: $r|_W<0$, $r|_{\partial W}=0$,
$dr\ne 0$ near $\partial W$. Let
$j:\partial W \hookrightarrow \overline{W}$ be the inclusion map.
The $1$-form $\alpha=j^*{\rm Im}(\overline{\partial r})$
is a contact form on~$\partial W$.
Denote by $\nu$ the measure on $\partial W$ associated
to the $(2n-1)$-form $\alpha\wedge (d\alpha)^{n-1}$,
and denote $L^2=L^2(\partial W,\nu)$.
Denote by $A(W)$ the space of functions on~$\overline{W}$
which are continuous on~$\overline{W}$, smooth
on $\partial W$, and holomorphic on~$W$.
Def\/ine {\it the Hardy space} $H^2=H^2(\partial W)$
to be the closure in $L^2$ of $\{ f|_{\partial W} \ |
\ f\in A(W)\}$.
Denote by $\Pi:L^2\to H^2$ the orthogonal projector.

By def\/inition an operator $T:C^{\infty}(\partial W)\to
C^{\infty}(\partial W)$ is called a {\it Toeplitz
operator of order $k$} if it is of the form
$\Pi Q\Pi$, where $Q$ is a pseudodif\/ferential operator
of order $k$. The {\it symbol} of $T$
is $\sigma(T):=\sigma(Q)|_{\Sigma}$
(a function on $\Sigma$), where
$\sigma (Q)$ is the symbol of $Q$ and
\[
\Sigma=\{ (x,\xi)\, |\, x\in\partial W,\; \xi=r\alpha _x,\; r>0\}
\]
is a symplectic submanifold of $T^*\partial W$.
Note that the symbol is well-def\/ined:
if $Q_1$, $Q_2$ are pseudodif\/ferential operators and
$T=\Pi Q_1\Pi=\Pi Q_2\Pi$, then
$\sigma(Q_1)|_{\Sigma}=\sigma(Q_2)|_{\Sigma}$.
Boutet de Monvel and Guillemin also show that
for Toeplitz operators $T_1$, $T_2$
$\sigma(T_1T_2)=\sigma(T_1)\sigma(T_2)$,
$\sigma([T_1,T_2])=\{ \sigma(T_1),\sigma(T_2)\}$,
where $\{\cdot,\cdot\}$ is the intrinsic Poisson bracket
on the symplectic mani\-fold $\Sigma$, and
Toeplitz operators form a ring under composition.

\subsection[Berezin-Toeplitz quantization]{Berezin--Toeplitz quantization}
\label{bt}

Let $X$ be a connected compact
$n$-dimensional K\"ahler manifold, $n\ge 1$.
Denote the K\"ahler form by $\omega$.
Assume that $\omega$ is integral.
There is an (ample) holomorphic Hermitian line bundle $L\to X$,
with Hermitian connection $\nabla$, such that
$curv (\nabla)=-2\pi i\omega$ (thus $c_1(L)=[\omega]$).
Let $N$ be a positive integer.
We shall denote by $L^2(X,L^{\otimes N})$ the space
of square-integrable sections of~$L^{\otimes N}$,
and by $H^0(X,L^{\otimes N})$ the space of holomorphic
sections of $L^{\otimes N}$.
Also denote by
$C^\infty (X)$ the space of
real-valued smooth functions on $X$ and by
$C_{{\mathbb C}}^\infty (X)$ the space of
complex-valued smooth functions on $X$.

The unit disc bundle $W$ in $L^*$ is a strictly
pseudoconvex domain. Denote $P=\partial W$
(the unit circle bundle in $L^*$).
In this particular setting, with $k=0$,
the def\/inition of a Toeplitz operator in Section
\ref{bg}
leads to the following (revised) def\/inition.
Let  $f\in C_{{\mathbb C}}^\infty (X)$.
The corresponding {\it Toeplitz operator}
(also called {\it Berezin--Toeplitz operator}
in this setting) is
\[
T_f=\oplus_{N=0}^{\infty}T_f^{(N)},
\]
where
$T_f^{(N)}=\Pi^{(N)}\circ M_f^{(N)}\in
{\rm End}(H^0(X,L^{\otimes N}))$,
\[
M_f^{(N)}: \ H^0(X,L^{\otimes N})\to L^2(X,L^{\otimes N}),\qquad
s\mapsto fs,
\]
is the operator of multiplication by $f$
and
\[
\Pi^{(N)}:\ L^2(X,L^{\otimes N})\to H^0(X,L^{\otimes N})
\]
is the orthogonal projector.

Much has been written on this subject, see, in particular,
review papers \cite{borthwick:00, schlichenmaier:95,schlichenmaier:01}.

We shall list some properties of Berezin--Toeplitz operators.
For every $N$ the map
\[
C_{{\mathbb C}}^\infty (X)\to {\rm End}\big(H^0(X,L^{\otimes N})\big)
,\qquad
f\mapsto T^{(N)}_f
\]
is surjective~\cite[Proposition 4.2]{bordemann:94}.

It is known that for a positive integer $m$,
$f_1,\dots,f_m\in C^{\infty}(X)$
\[
{\rm tr} \big(T_{f_1}^{(N)}\cdots T_{f_m}^{(N)}\big)=
N^n\left(\int_X f_1\cdots f_m\frac{\omega^n}{n!}+O\left(\frac{1}{N}\right)\right)
\]
as $N\to +\infty$ \cite[Section 5]{bordemann:94}.

Theorem 4.2 \cite{bordemann:94}
states that for $f,g\in C^\infty (X)$
\[
\big\|N[T^{(N)}_f,T^{(N)}_g]-iT^{(N)}_{\{ f,g\} }\big\|=O\left(\frac{1}{N}\right)
\]
as $N\to +\infty$, where
$\|\cdot \|$ is the operator norm, i.e.
$||A^{(N)}||^2=\sup _{s\in(H^0(X,L^{\otimes N})-0)}
\frac{\langle As,As\rangle_N}{\langle s,s\rangle_N}$
for $A^{(N)}\in {\rm End}(H^0(X,L^{\otimes N}))$, and
$\langle \cdot,\cdot\rangle _N$ is the Hermitian inner product
on $H^0(X,L^{\otimes N})$.
Simi\-lar statements, for certain deformations
of Lie algebra structure on ${\rm End}(H^0(X,L^{\otimes N}))$,
were obtained in~\cite{foth:07}.

Also
\[
\big\|T_{fg}^{(N)}-T_f^{(N)}T_g^{(N)}\big\|=O\left(\frac{1}{N}\right)
\]
 as $N\to +\infty$
\cite[p.~291, (2)]{bordemann:94}.

\subsection{Other aspects}

{\bf 2.3.1.}
There is a strong connection between Berezin--Toeplitz
quantization and deformation quantization.
See, in particular,
\cite{borthwick:93,charles:03,guillemin:95,karabegov:01,
klimek:92,klimek:96,reshetikhin:00,schlichenmaier:00}.

\medskip

\noindent
{\bf 2.3.2.}
Everything discussed in Section~\ref{bt}
has a traditional ``translation'' into the
language of physics: it is customary to say
that $X$ is the classical phase space,
$1/N$ is the Planck's constant, $H^0(X,L^{\otimes N})$
is the space of wave functions,
$f$ is a classical Hamiltonian,
$T_f^{(N)}$ is the quantum Hamiltonian,
and $N\to +\infty$ is the semiclassical limit.
There are interesting and dif\/f\/icult results
related to semiclassical behaviour of eigenvalues
and eigenfunctions of Toeplitz operators,
to quantization of maps,
and, generally, to the relation between
classical dynamics and quantum dynamics.
See, in particular,
\cite{borthwick:98,boutet:85,zelditch:97,
zelditch:05,zelditch:95}.

\medskip

\noindent
{\bf 2.3.3.}
Symbol calculus of Toeplitz operators
has been used to study integrable systems
\cite{bloch:03}.

\medskip

\noindent
{\bf 2.3.4.}
Lagrangian submanifolds and symplectic reduction
are two very important concepts in~symplectic geometry.
Toeplitz operators in this context have been
studied and exploited  in~\cite{borthwick:95,charlesl:03,charles:06,
charlesl:06,charles:07}.

\medskip

\noindent
{\bf 2.3.5.}
Let $(X,\omega,J_0)$ be a compact connected K\"ahler manifold
(here $X$ denotes the underlying smooth manifold,
$\omega$ is the symplectic structure, $J_0$ is the complex
structure). Denote by ${\mathcal J}$ the space
of all complex structures on $X$ compatible with $\omega$,
it is an inf\/inite-dimensional K\"ahler manifold.
The group ${\rm Symp}(X,\omega)$ acts on ${\mathcal J}$
by pull-back.
It was observed by Fujiki \cite{fujiki:92} and
by Donaldson \cite{donaldson:97}
that $J\to s_J$, where $J\in {\mathcal J}$
and $s_J$ is the scalar curvature of the
Riemannian metric given by $\omega$ and $J$,
is a moment map for this action.
In \cite{foth:08} we obtain another proof
of this statement, using Toeplitz operators
and semiclassical asymptotics.

\subsection*{Acknowledgements}
This work supported in part by NSERC.

\pdfbookmark[1]{References}{ref}
\LastPageEnding


\begin{thebibliography}{99}

\footnotesize\itemsep=0pt

\bibitem{berezin:72}Berezin F.,
Covariant and contravariant symbols of operators,
{\it Izv. Akad. Nauk SSSR Ser. Mat.} {\bf 36} (1972), 1134--1167 (in Russian).

\bibitem{berezin:74}Berezin F.,
Spectral properties of generalized Toeplitz matrices,
{\it Mat. Sb. (N.S.)} {\bf 95 (137)} (1974), 305--325, 328  (in Russian).

\bibitem{berezin:75}Berezin F.,
General concept of quantization,
{\it Comm. Math. Phys.} {\bf 40} (1975), 153--174.

\bibitem{berman:05}Berman R.,  Berndtsson B.,  Sjoestrand J.,
Asymptotics of Bergman kernels,
\href{http://arxiv.org/abs/math.CV/0506367}{math.CV/0506367}.

\bibitem{bloch:03}Bloch A.,  Golse F.,  Uribe A.,
Dispersionless Toda and Toeplitz operators,
{\it Duke Math. J.} {\bf 117} (2003), 157--196.

\bibitem{bordemann:94}Bordemann M.,  Meinrenken E.,
 Schlichenmaier M.,
Toeplitz quantization of K\"ahler manifolds and $gl(N)$, $N\to\infty$ limits,
{\it Comm. Math. Phys.} {\bf 165} (1994), 281--296.

\bibitem{borthwick:00}Borthwick D.,
Introduction to K\"ahler quantization, in Quantization,
the Segal--Bargmann transform and semiclassical analysis,
1st Summer School in Analysis and Mathematical Physics (Mexico, 1998),
{\it Contemp. Math.} {\bf 260} (2000),
91--132.

\bibitem{borthwick:93}Borthwick D.,  Lesniewski A.,
 Upmeier H.,
 Nonperturbative deformation quantization of Cartan domains,
{\it J. Funct. Anal.} {\bf 113} (1993), 153--176.

\bibitem{borthwick:98}Borthwick D., Paul T., Uribe A.,
Semiclassical spectral estimates for Toeplitz operators,
{\it Ann. Inst. Fourier (Grenoble)} {\bf 48} (1998), 1189--1229.

\bibitem{borthwick:95}Borthwick D.,  Paul T.,  Uribe A.,
 Legendrian distributions with applications to relative
Poincar\'e series,
{\it Invent. Math.} {\bf 122} (1995), 359--402, \href{http://arxiv.org/abs/hep-th/9406036}{hep-th/9406036}.

\bibitem{bottcher:90}B\"ottcher A.,  Silbermann B.,
Analysis of Toeplitz operators,
Springer-Verlag, Berlin, 1990.

\bibitem{boutet:78}Boutet de Monvel L.,
On the index of Toeplitz operators of several complex variables,
{\it Invent. Math.} {\bf 50} (1978/79), 249--272.

\bibitem{boutet:85}Boutet de Monvel L.,
Toeplitz operators~-- an asymptotic quantization of symplectic cones, in
Stochastic Processes and Their Applications in Mathematics and Physics (Bielefeld, 1985),
{\it  Math. Appl.},  Vol.~61,
Kluwer Acad. Publ., Dordrecht, 1990, 95--106.

\bibitem{boutet:05}Boutet de Monvel L.,
Logarithmic trace of Toeplitz projectors,
{\it Math. Res. Lett.} {\bf 12} (2005),  401--412, \href{http://arxiv.org/abs/math.CV/0412252}{math.CV/0412252}.

\bibitem{boutet:81}Boutet de Monvel L.,  Guillemin V.,
The spectral theory of Toeplitz operators,
{\it Annals of Math. Studies}, Vol.~99,
Princeton University Press, Princeton,
New Jersey, 1981.

\bibitem{boutet:76}Boutet de Monvel L.,  Sj\"ostrand J.,
Sur la singularite des noyaux de Bergman et de Szego,
in Journ\'ees: \'Equations aux D\'erivees Partielles de Rennes (1975),
{\it Asterisque}, no.~34--35 (1976), 123--164. 

\bibitem{brown:63}Brown A.,  Halmos P.,
 Algebraic properties of Toeplitz operators,
{\it J. Reine Angew. Math.}  {\bf 213} (1963/1964), 89--102.

\bibitem{charles:03}Charles L.,
 Berezin--Toeplitz operators, a semi-classical approach,
{\it Comm. Math. Phys.} {\bf 239} (2003),  1--28.

\bibitem{charlesl:03}Charles L.,
 Quasimodes and Bohr--Sommerfeld conditions for the Toeplitz operators,
{\it Comm. Partial Differential Equations} {\bf 28} (2003),
 1527--1566.

\bibitem{charles:06}Charles L.,
 Toeplitz operators and Hamiltonian torus actions,
{\it J. Funct. Anal.} {\bf 236} (2006), 299--350, \href{http://arxiv.org/abs/math.SG/0405128}{math.SG/0405128}.

\bibitem{charlesl:06}Charles L.,
 Symbolic calculus for Toeplitz operators with half-form,
{\it J. Symplectic Geom.} {\bf 4} (2006), 171--198, \href{http://arxiv.org/abs/math.SG/0602167}{math.SG/0602167}.

\bibitem{charles:07}Charles L.,
 Semi-classical properties of geometric quantization with metaplectic correction,
{\it Comm. Math. Phys.}  {\bf 270} (2007),  445--480, \href{http://arxiv.org/abs/math.SG/0602168}{math.SG/0602168}.

\bibitem{donaldson:97}Donaldson S.,
 Remarks on gauge theory, complex geometry and $4$-manifold topology, in Fields Medallists' Lectures,
{\it World Sci. Ser. 20th Century Math.}, Vol.~5, World Sci. Publ., River Edge, NJ, 1997, 384--403.

\bibitem{douglas:05}Douglas R.,
 Ideals in Toeplitz algebras,
{\it Houston J. Math.} {\bf 31} (2005), 529--539, \href{http://arxiv.org/abs/math.OA/0309438}{math.OA/0309438}.

\bibitem{douglas:70}Douglas R.,  Widom H.,
 Toeplitz operators with locally sectorial symbols,
{\it Indiana Univ. Math. J.} {\bf 20} (1970/1971), 385--388.

\bibitem{foth:07}Foth T.,
 Toeplitz operators, deformations, and asymptotics,
{\it J. Geom. Phys.} {\bf 57} (2007), 855--861.

\bibitem{foth:08}Foth T.,  Uribe A.,
The manifold of compatible almost complex structures and geometric quantization,
{\it Comm. Math. Phys.}, to appear.

\bibitem{fujiki:92}Fujiki A.,
 Moduli space of polarized algebraic manifolds and K\"ahler metrics,
{\it Sugaku Expositions} {\bf 5} (1992), 173--191.

\bibitem{guillemin:95}Guillemin V.,
 Star products on compact pre-quantizable symplectic manifolds,
{\it Lett. Math. Phys.} {\bf 35} (1995), 85--89.

\bibitem{jewell:80}Jewell N.,
 Toeplitz operators on the Bergman spaces
and in several complex variables,
{\it Proc. London Math. Soc. (3)} {\bf 41} (1980), 193--216.

\bibitem{jewell:79}Jewell N.,  Krantz S.,
Toeplitz operators and related
function algebras on certain pseudoconvex domains,
{\it Trans. Amer. Math. Soc.} {\bf 252} (1979), 297--312.

\bibitem{karabegov:01}Karabegov A.,  Schlichenmaier M.,  Identif\/ication of Berezin--Toeplitz deformation quantization,
{\it J. Reine Angew. Math.} {\bf 540} (2001), 49--76, \href{http://arxiv.org/abs/math.QA/0006063}{math.QA/0006063}.

\bibitem{klein:72}Klein E.,
The numerical range of a Toeplitz operator,
{\it Proc. Amer. Math. Soc.} {\bf 35} (1972), 101--103.

\bibitem{klein:73}Klein E.,
 More algebraic properties of Toeplitz operators,
{\it Math. Ann.} {\bf 202} (1973), 203--207.

\bibitem{klimek:92}Klimek S.,  Lesniewski A.,
 Quantum Riemann surfaces. I.~The unit disc,
{\it Comm. Math. Phys.} {\bf 146} (1992), 103--122.\\
Klimek S.,  Lesniewski A.,
Quantum Riemann surfaces. II.~The discrete series,
{\it Lett. Math. Phys.} {\bf 24} (1992),  125--139.\\
Klimek S.,  Lesniewski A., Quantum Riemann surfaces. III.~The exceptional cases,
{\it Lett. Math. Phys.} {\bf 32} (1994), 45--61.

\bibitem{klimek:96}Klimek S.,  Lesniewski A.,
 Quantum Riemann surfaces for arbitrary Planck's constant,
{\it J. Math. Phys.} {\bf 37} (1996), 2157--2165.

\bibitem{kostant:70}Kostant B.,
 Quantization and unitary representations. I.~Prequantization,
in Lectures in Modern Analysis and Applications.~III,  {\it Lecture Notes in Math.}, Vol.~170,
Springer, Berlin, 1970, 87--208.

\bibitem{peller:93}Peller V.,
 Invariant subspaces of Toeplitz operators with piecewise continuous symbols,
{\it Proc. Amer. Math. Soc.} {\bf 119} (1993), 171--178.


\bibitem{reshetikhin:00}Reshetikhin N.,  Takhtajan L.,
 Deformation quantization of K\"ahler manifolds,
in L.D.~Faddeev's Seminar on Mathematical Physics,,
{\it Amer. Math. Soc. Transl. Ser. 2}, Vol.~201,
 Amer. Math. Soc., Providence, RI, 2000, 257--276, \href{http://arxiv.org/abs/math.QA/9907171}{math.QA/9907171}.

\bibitem{salinas:89}Salinas N.,  Sheu A.,  Upmeier H.,
Toeplitz operators on pseudoconvex domains and foliation
$C^*$-algebras,
{\it Ann. Math. (2)} {\bf 130} (1989),  531--565.

\bibitem{schlichenmaier:95}Schlichenmaier M.,
 Berezin--Toeplitz quantization of compact K\"ahler manifolds,
in Quantization, Coherent States, and Poisson Structures (Bialowieza, 1995),
PWN, Warsaw, 1998, 101--115.

\bibitem{schlichenmaier:00}Schlichenmaier M.,
 Deformation quantization of compact K\"ahler manifolds by Berezin--Toeplitz quantization, in
 Conf\'erence Mosh\'e Flato 1999, Vol.~II (Dijon),
{\it Math. Phys. Stud.}, Vol.~22, Kluwer Acad. Publ.,
Dordrecht, 2000, 289--306, \href{http://arxiv.org/abs/math.QA/9910137}{math.QA/9910137}.

\bibitem{schlichenmaier:01}Schlichenmaier M.,
 Berezin--Toeplitz quantization and Berezin transform,
in Long Time Behaviour of Classical and Quantum Systems
(Bologna, 1999),  {\it Ser. Concr. Appl. Math.},
Vol.~1, World Sci. Publ., River Edge, NJ, 2001, 271--287.

\bibitem{souriau:97}Souriau J.-M.,
Structure of dynamical systems,
A symplectic view of physics,
{\it Progress in Mathematics}, Vol.~149,
Birkh\"auser Boston, Inc., Boston, MA, 1997 (translation of  Structure des syst\`emes dynamiques,
Ma\^itrises de math\'ematiques Dunod, Paris, 1970).

\bibitem{stroethoff:03}Stroethof\/f K.,  Zheng D.,
Bounded Toeplitz products on the Bergman space of the polydisk,
{\it J. Math. Anal. Appl.} {\bf 278} (2003), 125--135.

\bibitem{stroethoff:07}Stroethof\/f K.,  Zheng D.,
Bounded Toeplitz products on Bergman spaces of the unit ball,
{\it J. Math. Anal. Appl.} {\bf 325} (2007), 114--129.

\bibitem{upmeier:83}Upmeier H.,
 Toeplitz operators on bounded symmetric domains,
{\it Trans. Amer. Math. Soc.} {\bf 280} (1983), 221--237.

\bibitem{upmeier:84}Upmeier H.,
 Toeplitz $C^{*} $-algebras on bounded symmetric domains,
{\it Ann. Math. (2)}  {\bf 119} (1984),  549--576.

\bibitem{upmeier:85}Upmeier H.,
 Toeplitz operators on symmetric Siegel domains,
{\it  Math. Ann.} {\bf 271} (1985), 401--414.

\bibitem{upmeier:96}Upmeier H.,
Toeplitz operators and index theory in several complex variables,
{\it Operator Theory: Advances and Applications}, Vol.~81,
 Birkh\"auser Verlag, Basel, 1996.

\bibitem{widom:64}Widom H.,
 On the spectrum of a Toeplitz operator,
{\it Pacific J. Math.} {\bf 14} (1964), 365--375.

\bibitem{widom:66}Widom H.,
Toeplitz operators on $H_{p}$,
{\it Pacific J. Math.} {\bf 19} (1966), 573--582.

\bibitem{zelditch:97}Zelditch S.,
 Index and dynamics of quantized contact transformations,
{\it Ann. Inst. Fourier (Grenoble)}
{\bf 47} (1997), 305--363, \href{http://arxiv.org/abs/math-ph/0002007}{math-ph/0002007}.

\bibitem{zelditch:05}Zelditch S.,
 Quantum maps and automorphisms,
in The Breadth of Symplectic and Poisson Geometry,
{\it Progr. Math.}, Vol.~232,
Birkh\"auser Boston, Boston, MA, 2005, 623--654, \href{http://arxiv.org/abs/math.QA/0307175}{math.QA/0307175}.

\bibitem{zelditch:95}Zelditch S.,
Quantum dynamics from the semi-classical
point of view, unpublished notes, available at\linebreak \url{http://mathnt.mat.jhu.edu/zelditch/}.





\end{thebibliography}
\end{document}